\documentclass[letterpaper, conference]{IEEEtran}
\IEEEoverridecommandlockouts

\usepackage{algorithm}
\usepackage[noend]{algpseudocode}

\usepackage[noadjust]{cite}

\usepackage{amsmath}
\usepackage{amssymb}
\usepackage{amsthm}
\usepackage{amsfonts}
\usepackage{nccmath}

\usepackage[noabbrev]{cleveref}
\usepackage{nameref}

\newtheorem{theorem}{Theorem}

\newtheorem{lemma}[theorem]{Lemma}

\theoremstyle{definition}

\newtheorem{assumption}[theorem]{Assumption}

\theoremstyle{remark}
\newtheorem{remark}[theorem]{Remark}

\usepackage{mathtools}
\usepackage{graphicx}
\graphicspath{ {./images/} }
\usepackage{subcaption}
\usepackage{color}
\usepackage{xcolor}
\usepackage{float}

\usepackage{balance}
\crefformat{equation}{(#1)}
\usepackage{bbm}
\usepackage{bm}
\usepackage{comment}
\usepackage{arcs}

\usepackage{resizegather}

\usepackage{ifthen}
\newcommand{\prob}[2][]{
\ifthenelse { \equal {#2} {} }  
    { \mathbb{P}_{#1}}   
    { \mathbb{P}_{#1}\left[#2\right] }   
}

\usepackage{acronym}
\newacro{lqr}[LQR]{Linear Quadratic Regulators}
\newacro{slqr}[SLQR]{Structured Linear Quadratic Regulators}
\newacro{olqr}[OLQR]{Output-feedback Linear Quadratic Regulators}
\newacro{lqg}[LQG]{Linear Quadratic Gaussian}
\newacro{dare}[DARE]{Discrete-time Algebraic Riccati Equation}
\newacro{ouralgo}[RNPO]{Riemannian Newton-type Policy Optimization}
\newacro{PO}[PO]{Policy Optimization}
\newacro{pg}[PG]{Projected Gradient}
\newacro{mdp}[MDP]{Markov Decision Process}

\usepackage{textcomp}
\def\BibTeX{{\rm B\kern-.05em{\sc i\kern-.025em b}\kern-.08em
		T\kern-.1667em\lower.7ex\hbox{E}\kern-.125emX}}

\usepackage{enumerate}

\usepackage{paralist}

\def\x{{X}}
\def\s{{\Lambda}}

\def\w{{W}}
\def\u{{U}}
\def\H{{\mathcal{H}}}

\def\E{{\mathbb{E}\,}}

\def\P{\mathbb{P}}

\def\F{\mathcal{F}}

\def\R{\mathbb{R}}

\def\Amatrices{\mathbb{R}^{n\times n}}
\def\Bmatrices{\mathbb{R}^{n\times m}}

\def\Hmatrices{\mathbb{R}^{n\times d}}

\def\stableK{\mathcal{S}}

\newcommand{\tr}[1]{\ensuremath{\mathrm{tr}\left[ #1 \right]}}

\usepackage{tikz}
\allowdisplaybreaks

\title{\vspace{0.5cm}Ergodic-Risk Constrained Policy Optimization:\\ The Linear Quadratic Case
\thanks{The authors are with the School of Engineering and Applied Sciences, Harvard University, Cambridge, USA. Emails: \textit{talebi@seas.harvard.edu} and  \textit{nali@seas.harvard.edu}. This work is supported by NSF AI Institute 2112085. }
}
\author{ Shahriar Talebi,  \ Na Li

}

\begin{document}
\maketitle

\begin{abstract}
Risk-sensitive control balances performance with resilience to unlikely events in uncertain systems. This paper introduces \textit{ergodic-risk criteria}, which capture long-term cumulative risks through probabilistic limit theorems.
By ensuring the dynamics exhibit strong ergodicity, we demonstrate that the time-correlated terms in these limiting criteria converge even with potentially heavy-tailed process noises as long as the noise has a finite fourth moment. Building upon this, we proposed the ergodic-risk constrained policy optimization which incorporates an ergodic-risk constraint to the classical Linear Quadratic Regulation (LQR) framework. We then propose a primal-dual policy optimization method that optimizes the average performance while satisfying the ergodic-risk constraints. Numerical results demonstrate that the new risk-constrained LQR not only optimizes average performance but also limits the asymptotic variance associated with the ergodic-risk criterion, making the closed-loop system more robust against sporadic large fluctuations in process noise.

\end{abstract}

\begin{IEEEkeywords}
	    \textit{Ergodic-risk; Risk-aware Optimal Control; Risk-averse Decision Making; Linear Quadratic Regulator (LQR); Uniformly Ergodic Chains; Constrained Policy Optimization}
\end{IEEEkeywords}


\section{Introduction}
\label{sec:intro}

Optimizing average performance, as is typical in standard stochastic optimal control, often fails to yield effective policies for decision making in stochastic environments where deviations from expected outcomes carry significant risk; e.g. in financial markets \cite{rockafellar_optimization_2000}, safe robotics and autonomous systems \cite{majumdar_how_2020}, and healthcare \cite{eichler_risks_2013}.
As such, incorporating risk measures become vital in such decision making problems for balancing the performance with resilience to rare events.

While robust control frameworks (e.g. the mixed $\mathcal{H}_2$-$\mathcal{H}_\infty$ in \cite{zhang_policy_2021}), focus on incorporating the worst-case scenario performance (e.g. $\mathcal{H}_\infty$-norm) as a constraint, they can be overly conservative (or occasionally unfeasible) when those unlikely events are (possibly) unbounded. Risk-aware approaches, on the other hand, offer a (probabilistic) compromise by building on available stochastic priors to manage both risk and performance, simultaneously. 
Consequently, there have been significant efforts \cite{whittle_risk-sensitive_1981,borkar_risk-constrained_2014,chow_risk-constrained_2018,sopasakis_risk-averse_2019} in developing risk-aware decision making frameworks using tools like Conditional Value at Risk (CVaR) \cite{rockafellar_optimization_2000}, Markov Risk Measure \cite{ruszczynski_risk-averse_2010}, and Entropic Value at Risk (EVaR) \cite{ahmadi-javid_entropic_2012}, offering a better balance by considering both risk and average performance. 

These measures are often deployed on finite-horizon variables \cite{sopasakis_risk-averse_2019}, under finite first-hitting time \cite{chow_risk-constrained_2018}, and/or in finite state-space models \cite{borkar_risk-constrained_2014}. While these settings avoid complications regarding limiting probabilities, they may not fully capture \emph{long-term risk} associated with the stochastic behaviors, especially in \emph{unbounded general-state} Markov processes--see the recent survey \cite{biswas_ergodic_2023}.
Also at the stationary limits of the process, the optimal policy that minimizes the worst-case CVaR of the quadratic cost is shown to be equivalent to that of the \ac{lqr} optimal \cite{kishida_risk-aware_2023}. This is yet another evidence suggesting that risk-sensitive design is particularly critical in \emph{nonstationary} processes, where their statistical properties are still changing over time.
Among these, the folklore risk-sensitive framework by Whittle \cite{whittle_risk-sensitive_1981}, aka Linear Exponential Quadratic Gaussian (LEQG), handles the general unbounded, nonstationary setting--which (in certain parameter regimes) can be interpreted as optimizing a specific mixture of the average performance and its higher moments (e.g. variance). However, the Gaussian noise (with finite moments of all orders) is critical for the exponentiation to be well-defined, and thus does not capture cases with \textit{heavy-tailed noise} distributions modeling rare events. This motivated \cite{tsiamis_risk-constrained_2020} to introduce a framework for constraining the uncertainty in the ``state-related portion'' of the finite-horizon \ac{lqr} cost, which is then extended to infinite-horizon through policy optimization techniques \cite{zhao_global_2023}.

This paper considers the stochastic Linear Time Invariant (LTI) model with (unbounded) \textit{heavy-tailed} process noise, providing a framework for risk-aware decision-making in \textit{unbounded}, \textit{nonstationary} Markov processes.
%
We introduce \emph{ergodic-risk criteria} to address risks in the long-term stochastic behavior, accounting for extreme deviations beyond mean performance (\Cref{sec:probSetup}). 
Built upon this, we propose the ergodic-risk constrained policy optimization which incorporates the ergodic-risk constraints in the classical LQR framework. 
By ensuring a strong ergodicity of the process \cite{meyn_markov_2009}, we handle system's state correlations and characterize \textit{quadratic} ergodic-risk criteria as long as the process noise has finite fourth moment (\Cref{sec:F-CLT}). 
This enables us to address long-term risk in non-stationary processes, previously excluded from the literature (see e.g. \cite{kishida_risk-aware_2023}).
We then consider a primal-dual policy optimization method based on strong duality that optimizes the average performance while meeting the risk constraints (\Cref{sec:primal-dual}), with convergence guarantee in \Cref{subsec:convergence}. 
Finally, we demonstrate the numerical performance of the algorithm over randomly sampled problem instances and contrast the quality of the ergodic-risk optimal policy against the LQR optimal policy (\Cref{subsec:simulations}). 

Finally, our ergodic-risk analysis is fundamentally different than that of \cite{tsiamis_risk-constrained_2020,zhao_global_2023} and in fact, the proposed Ergodic-risk framework (\Cref{thm:C-infty-N-convergence}) generalizes to provably capture the long-term risk associated with heavy-tailed noise for any quadratic functional of \emph{both} the states and the inputs. However similar to \cite{tsiamis_risk-constrained_2020,zhao_global_2023}, the strong duality (\Cref{sec:strong-duality}) and the resulting primal-dual algorithm (\Cref{sec:primal-dual-algo}) are established only for the case in which the risk functional does not explicitly depend on the input signal. 

\section{Ergodic-Risk Criteria and Problem Setup}
\label{sec:probSetup}
Consider the discrete-time stochastic linear system,
\begin{equation}\label{eq:dynamics}
    \x_{t+1} = A \x_{t} + B \u_{t} + H \w_{t+1}, \quad t\geq 0,
\end{equation}
where $A \in \Amatrices$, $B \in \Bmatrices$, and $H \in \Hmatrices$ are system parameter matrices, $\x_t$ and $\u_t$ denote the stochastic state and input (vectors), respectively. Also, $\w_t$ is denoting the process noise and $\x_0$ is the initial state vectors. $\w_t$ and $\x_0$ are independently sampled from zero-mean probability distributions $\prob[\w]{}$ and $\prob[0]{}$ with covariances $\Sigma_W$ and $\Sigma_0$, respectively. The history of state-input trajectory up to time $t$ is denoted by $\H_t = \{\x_j,\u_j\}_{j=0}^{t}$, and $\F_t = \sigma(\H_t)$ denotes the $\sigma$-algebra generated by $\H_t$. 
Each $\w_t$ is measurable with respect to $\F_{t}$, denoted by $\w_t \in\F_{t}$.
Let $\F_{-1}$ denote the trivial $\sigma$-algebra, and $\F_{t-1}\vee \sigma\{\x_t\}$ denote the smallest $\sigma$-algebra containing both $\F_{t-1}$ and $\sigma\{\x_t\}$. At each time $t\geq 0$, we can apply an \textit{admissible} input $\u_t \in \mathcal{L}_2(\F_{t-1}\vee \sigma\{\x_t\})$, i.e. a square-integrable, measurable function with respect to $\F_{t-1}\vee \sigma\{\x_t\}$, and then measure the next state $\x_{t+1}$. 
We restrict ourselves to \textit{stationary Markov policies} $\pi$, measurable mappings independent of time $t$ that generate an input sequence $\{\u_t = \pi(\H_{t-1},\x_t)\}_t$ such that $\u_t$ is admissible for all time $t$ and does not depend on the history $\H_{t-1}$--see \cite{meyn_markov_2009,durrett_probability_2019} for further details regarding probabilistic notions.

We require the process noise and initial states to be uncorrelated across time, i.e.,
$\E[\x_0 \w_t^\intercal]= 0$, and $\E[\w_{t+\tau} \w_t^\intercal]=0, \; \forall t,\tau \geq 1$, so that $\w_t$ is independent of $\F_{t-1}$ for all $t\geq1$.
For simplicity, we pose the following assumption:
\begin{assumption}\label{assmp:noise}
    The sequence $\{\w_t\}$ is i.i.d. samples of a common zero-mean probability measure $\P_W$ that is non-singular with respect to Lebesgue measure on $\R^d$ and has a non-trivial density, with a finite covariance $\Sigma_W \succ 0$.
    \qed
\end{assumption}

Given an (admissible) input signal $\u = \{\u_t\}$, we define the cumulative \emph{performance cost}
\(
    J_{T}(\u) = \sum_{t=0}^T \x_t^\intercal Q \x_t +  \u_t^\intercal R \u_t,
\)
with $Q \succeq 0$ and
$R\succ 0$ being positive semidefinite and positive definite matrices,
respectively. The standard infinite-horizon \ac{lqr} problem is then to design a sequence of admissible inputs $\u=\{\u_t\}_0^\infty$ that minimizes
\begin{equation}\label{eq:costx0}
\textstyle
     J(\u) =  \limsup_{T\to\infty} \frac{1}{T} \E [J_T(\u)],
\end{equation}
subject to dynamics in \Cref{eq:dynamics}.
It is well known \cite{goodwin_control_2001} 
that the optimal solution reduces to solving the \ac{dare} for a cost matrix $P_{\mathrm{LQR}}$ and the LQR optimal input is $\u_t^*= K_{\mathrm{LQR}} \x_t$ with the controller (or ``policy'') 
\(K_{\mathrm{LQR}} = - (R + B^\intercal P_{\mathrm{LQR}} B)^{-1} B^\intercal P_{\mathrm{LQR}}  A.\)
So, the optimal LQR policy is a \emph{linear} stationary Markov policy, and we restrict ourselves to the same class in this work.

\subsection{Ergodic-Risk Criterion for Stochastic Systems}
We propose a risk criterion that captures the long-term accumulative uncertainty by adding step-wise uncertainties as the system evolves. Since each state $\x_t$ is observed iteratively, it is natural to consider the uncertainty at each stage and accumulate these contributions over time to characterize the overall risk. This leads to a cumulative uncertainty variable, which converges to a limiting value if properly normalized.
 
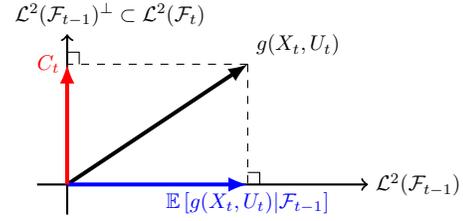
\begin{figure}[!pt]
    \centering
    \begin{tikzpicture}[scale=0.8, transform shape] 
    \draw[thick,->] (-0.5,0) -- (5,0) node[right] {$\mathcal{L}^2(\mathcal{F}_{t-1})$};
    \draw[thick,->] (0,-0.5) -- (0,2.5) node[above right] {\hspace{-1cm}$\mathcal{L}^2(\mathcal{F}_{t-1})^\perp \subset \mathcal{L}^2(\mathcal{F}_{t})$};
    
    \coordinate (O) at (0,0); 
    \coordinate (V) at (3,2); 
    
    \draw[thick,->,line width=1.5pt,>=latex] (O) -- (V) node[above right] {$g(\x_{t},\u_{t})$};
    
    \draw[dashed] (V) -- (3,0) node[below] {\color{blue}$\E[g(\x_{t},\u_{t}) | \F_{t-1}]$};
    \draw[dashed] (V) -- (0,2) node[left] {\color{red}$C_t$};
    
    \draw (3.2,0) -- ++(0,0.2) -- ++(-0.2,0); 
    \draw (0,2.2) -- ++(0.2,0) -- ++(0, -0.2); 
    
    \draw[thick,<-,blue,line width=1.5pt,>=latex] (3,0) -- (0,0);  
    \draw[thick,<-,red,line width=1.5pt,>=latex] (0,2) -- (0,0);  
\end{tikzpicture}
    \caption{\small The conditional expectation in blue is the orthogonal projection of $g(\x_{t},\u_{t})$ onto $\mathcal{L}^2(\mathcal{F}_{t-1})$, i.e. its best estimate by the information up to time $t-1$, solving $\arg\min_{\hat g \in \mathcal{L}^2(\mathcal{F}_{t-1})} \sqrt{\E[(g - \hat g)^2]}$. 
    So, $C_t$ (in red) then retains the ``uncertain component'' of $g(\x_{t},\u_{t})$.}
    \label{fig:one-step-risk}
    \vspace{-0.4cm}
\end{figure}

To formalize this, let us consider any measurable functional of choice $g:\R^n\times\R^m \mapsto \R$, called ``risk functional,'' for example a quadratic (or affine) function in $\x_t$ and $\u_t$ which evaluates the behavior of each sample path (possibly different than the performance cost $J_T$). At each time $t-1$, we have access to the past information in $\F_{t-1}$, so the risk factor at time $t$ is the ``uncertain component'' of the risk functional $g(\x_t,\u_t)$.
This motivates the following definition
\begin{equation}
    C_t \coloneqq 
    g(\x_{t},\u_{t}) -\E[g(\x_{t},\u_{t}) | \F_{t-1}], \quad \text{for $t\geq 0$,}
\end{equation}
 capturing the uncertainty in $g(\x_{t},\u_{t})$ at time $t$ relative to that past information---see \Cref{fig:one-step-risk}.
%
To account for the long-term risk behavior, especially in non-stationary processes with heavy-tailed noise, we define the \textit{ergodic-risk criterion}\footnote{A \emph{(uniformly) ergodic} Markov process visits all parts of the state space and uniformly converges to a unique stationary distribution, regardless of the starting point; see \cite{meyn_markov_2009}.} $C_\infty$ as the limit of the normalized cumulative uncertainty:
\begin{equation}\label{eq:def-C-infty}
\textstyle
    \frac{1}{\sqrt{t}} S_t \coloneqq \frac{1}{\sqrt{t}} \sum_{s=0}^t C_s \xrightarrow{d\;} C_\infty, \quad t\to\infty.
\end{equation}
We also consider the \emph{asymptotic conditional variance} $\gamma_N^2$ defined as the limit
\begin{equation}\label{eq:cond-var-def}
\textstyle
     \frac{1}{t} N_t\coloneqq \frac{1}{t} \sum_{s=1}^{t} \E[C_s^2|\F_{s-1}] \xrightarrow{a.s.} \gamma_N^2, \quad t\to \infty,
\end{equation}
which would serve as an ``estimate'' of the asymptotic variance of $C_\infty$, whenever well-defined (see \Cref{rem:asymp-var-conditional-var}).

\noindent \textbf{Problem.} We pose the so-called \textit{Ergodic-Risk Constrained Optimal Control Problem} (Ergodic-Risk COCP):
\begin{align}\label{eq:optimization}
   \displaystyle \min \; &J(U) \\
   \text{s.t.}~~ &\x_{t+1} = A \x_t + B U_t + H \w_{t+1}, \quad \forall t \geq0, \nonumber\\
  & \text{ constraints on risk measure over } C_\infty \nonumber.
\end{align}
As we will discuss further in \Cref{sec:primal-dual}, reasonable choices of constraints on (coherent) risk measure over $C_\infty$ or $\gamma_N^2$ (such as Conditional Value-at-Risk CVaR$_{\alpha}$ and the Entropic Value-at-Risk EVaR$_{\alpha}$ at a level $\alpha$), essentially reduces to an upperbound on a linear functional related to the variance of $C_\infty$, whenever is well-defined. Furthermore, as discussed later in \Cref{rem:asymp-var-conditional-var}, $\gamma_N^2$ will be used as the estimate of this variance, offering a more tractable form for policy optimization. 

The well-posedness of Ergodic-Risk COCP depends on the existence of $C_\infty$; so, we first discuss that $C_t$ has (at least) the \emph{necessary} properties to successfully capture this cumulative uncertainty.
Under continuity of risk functional $g$ and bounded moment conditions for the noise \cite[Lemma 2]{talebi_uniform_2024}, it can be shown that, $C_t$ is a Martingale Difference Sequence (MDS); i.e. for all $t$
    \begin{inparaenum}[(i)]
        \item $C_t$ is $\F_{t}$-adapted, and
        \item $\E |C_t| < \infty$, and
        \item $\E[C_t|\F_{t-1}] \stackrel{a.s.}{=} 0$;
    \end{inparaenum} however, guaranteeing that the limiting quantity $C_\infty$ is well-defined still requires careful analysis of convergence in distribution (and similarly for the limit in $\gamma_N^2$). 
Note that the summands in $C_\infty$ are highly correlated through the dynamics in \cref{eq:dynamics} and thus vanilla Central Limit Theorem does not apply as it requires independent summands. Even extended version of CLT for martingales \cite[Theorem 5.1]{komorowski_central_2012} is not directly useful here because the conditions translates to such strong stability conditions on \cref{eq:dynamics} that is not feasible by any feedback signal---unless the noise process $\{W_t\}$ is eventually vanishing, which is not the point of interest in this work.

Next, we study the ergodic-risk criteria defined in \cref{eq:def-C-infty} for quadratic risk functionals in LTI systems which will be used in solving the Quadratic Ergodic-Risk COCP in \Cref{sec:primal-dual}.

\section{Quadratic Ergodic-Risk for LTI systems}\label{sec:F-CLT}
Herein, we aim to quantify the ergodic-risk criteria $C_\infty$ and $\gamma_N^2$ defined in \cref{eq:def-C-infty} and \cref{eq:cond-var-def}. While we consider the \emph{quadratic} risk functionals $g$, the approach is more general. Also, the case of linear $g$ requires the noise to have only finite \emph{second} moment and follows similarly; thus is left to the reader.

For simplicity, we restrict ourselves to linear stationary Markov policies $\pi:x\mapsto Kx$ for some matrix parameter $K$, such that at each time $t$, the input is designed to be $\u_t = K \x_t$. We define the set of (Schur) stabilizing policies as 
\[\mathcal{S}=\left\{K \in \R^{m \times n}\;:\; A_K \coloneqq A+BK \text{ is Schur stable}\right\},\]
i.e. the closed-loop dynamics $A_K$ has spectral radius less than 1.
We refer to $\pi$ and $K$ as the policy without ambiguity. For $\mathcal{S}$ to be non-empty, we consider the minimal assumption:
\begin{assumption}\label{assmp:stability}
    The pair $(A,B)$ is stabilizable.
\end{assumption}
Let us define the following processes that become particularly relevant when the risk functional $g$ is quadratic;
namely, the running average and the running covariance

\begin{equation}\label{eq:def-s-t-Gamma-t}
\textstyle
     \s_t \coloneqq \sum_{s=1}^t \x_s , \text{ and }  \Gamma_t \coloneqq \sum_{s=1}^t \x_s \x_s^\intercal.
\end{equation}
Combining probabilistic tools (Markov's inequality) and system theoretic properties (exponential stability), one can characterize the expectation of running average and covariance, as well as boundedness in probability for the process $\{\x_t\}$ as summarized in the following result. Proof of this result combines \cite[Lemma 6]{talebi_data-driven_2023} with standard techniques and is deferred to the extended version \cite{talebi_uniform_2024}.

\begin{lemma}\label{lem:ave-conv}
    Under Assumptions \ref{assmp:noise} and \ref{assmp:stability}, for any stabilizing policy $K \in \mathcal{S}$, we have the following limits as $t\to\infty$:
    \(
        \E[\x_t] \to 0, \;\E[\s_t/t] \to 0, \text{ and } 
         \E[\Gamma_t/t] \to \Sigma_K,
    \)
    where $\Sigma_K$ is the unique positive definite solution to the Lyapunov equation:
    \begin{equation}\label{eq:Sigma-K}
        \Sigma_K = A_K \Sigma_K A_K^\intercal + H\Sigma_W H^\intercal.
    \end{equation}
    Furthermore, we obtain that
    \(\s_t/t \xrightarrow{p} 0, \text{ as } t\to\infty,\)
    and $\{\x_t\}$ is boundedness in probability. 
\end{lemma}

\Cref{lem:ave-conv} enables us to reason about the first and second moment of the process, however, it still does not provide enough for the convergence of $C_\infty$ in distribution. 
%
%
Here, we use another type of extensions to CLT known as ``Functional Central Limit Theorem'' that extends the Martingale CLT to Markov chains, connecting to their ``stochastic stability'' properties. It builds on the so-called ``uniform ergodicity'' \cite{meyn_markov_2009} as a stochastic stability notion that allows for such convergence to hold.
The next result shows convergence of $C_\infty$ for ``quadratic'' risk functionals which builds on
the Functional CLT recently tailored in \cite{talebi_uniform_2024}.

For that, let us consider the following stationary process: $\{Y_t: t\in\mathbb{Z}\}$ is the stationary process given by $Y_t = \sum_{n = 0}^\infty A_K^n H \w_{t-n}$ with $\{\w_t\}$ as i.i.d. samples of the same probability measure $\P_W$. Also, for any symmetric matrix $M = \sum_j \lambda_j v_j v_j^\intercal$, define
\(\textstyle \gamma_M^2 \coloneqq \sum_{j=1}^d \lambda_j^2 \gamma_{v_j}^2\)
where each 
\( \gamma_{v_j}^2 = \sum_{k=-\infty}^\infty \left( \E\left[(v_j^\intercal Y_0 Y_k^\intercal v_j)^2\right] - (v_j^\intercal \Sigma_K v_j)^2 \right).\)

\begin{theorem}\label{thm:C-infty-N-convergence}
    Suppose Assumptions \ref{assmp:noise} and \ref{assmp:stability} holds and consider the dynamics in \cref{eq:dynamics} for any policy $K \in \mathcal{S}$ that is stabilizing and $(A_K,H)$ is controllable. Consider
    \[g(x,u) = x^\intercal Q^c x + u^\intercal R^c u\]
    for some $Q^c, R^c\succeq 0$ and define $Q_K^c = Q^c + K^\intercal R^c K$. If the noise process $\{W_t\}$ has finite fourth moment, then the asymptotic variance $\gamma_M^2(K)$ with $M\coloneqq Q_K^c -A_K^\intercal Q_K^c A_K$ is well-defined, non-negative and finite.
    Then, if $\gamma_M^2(K)>0$,
    \[\textstyle \frac{1}{\sqrt{t}} S_t \coloneqq \frac{1}{\sqrt{t}} \sum_{s=1}^t C_s \xrightarrow{d\;} C_\infty \sim \mathcal{N}(0, \gamma_M^2(K)), \quad t\to\infty;\]
    otherwise, $\frac{1}{\sqrt{t}}\sum_{s=0}^t C_s \xrightarrow{a.s.} 0$.
    Furthermore,
    \[\textstyle \frac{1}{t} \sum_{s=1}^{t} C_s \xrightarrow{a.s.} 0,
    \quad
    \frac{1}{t} \sum_{s=1}^{t} \E [C_s^2| \F_{s-1}] \xrightarrow{a.s.} \gamma_N^2(K),\]
    where 
    \begin{gather*}
        \gamma_N^2(K) = 4\tr{Q_K^c H \Sigma_W H^\intercal Q_K^c (\Sigma_K-H \Sigma_W H^\intercal)}
        + m_4[Q_K^c],
    \end{gather*}
    with 
    \(
        m_4[Q_K^c] \coloneq \E\left[\tr{Q_K^c H (\w_{1} \w_{1}^\intercal - \Sigma_W )H^\intercal}^2\right].
    \)
\end{theorem}

\begin{remark}\label{rem:asymp-var-conditional-var}
First note that, by Tower property it also follows that
\(\frac{1}{t} \sum_{s=1}^{t} \E [C_s^2] \xrightarrow{a.s.} \gamma_N^2(K).\)
Second, the asymptotic conditional variance $\gamma_N^2(K)$ in \cref{eq:cond-var-def} serves as an ``estimate'' of $\gamma_M^2 = \lim_{t\to \infty} \E[S_t^2/t]$ with $S_t \coloneqq \sum_{s=1}^t C_s$ in the following sense: by Doob's decomposition we can show that
\(S_t^2 = M_t + N_t\)
where $(M_t,\F_t)$ is a martingale and $(N_t,\F_t)$ is predictable defined as
\(
    N_t - N_{t-1} = \E[S_t^2 - S_{t-1}^2|\F_{t-1}] 
    = \E[C_t^2|\F_{t-1}],
\)
where the second equality follows because $(S_t,\F_t)$ is a martingale.
Thus,
\(S_t^2  = M_t + \sum_{s=1}^{t} \E [C_s^2| \F_{s-1}].\)
So, $N_t$ can be interpreted as an intrinsic measure of time for the martingale $S_t$. Also, it can be interpreted as ``the amount of informaiton'' contained in the past history of the process, related to a standard Fisher information \cite[p. 54]{hall_martingale_1980}.
\end{remark}

\begin{proof}
    Consider the process $\{\x_t\}$ with any $K\in\mathcal{S}$ and starting from a fixed $\x_0 = x_0$.
    For simplicity, we define $G_{t} := g(\x_{t},\u_{t}), ~t\geq 0$ and note that $G_0 = x_0^\intercal Q_K^c x_0$ and
    \(
        G_{t+1} = (A_K\x_t+H\w_{t+1})^\intercal Q_K^c (A_K\x_t+H\w_{t+1}),
    \)
    for $t\geq 0$. 
     This, together with \Cref{assmp:noise} imply  that
    \[
        \E[G_{t+1}|\F_{t}] = \x_t^\intercal A_K^\intercal Q_K^c A_K\x_t + \tr{Q_K^c H \Sigma_W H^\intercal},
    \]
    where $\Sigma_W$ denotes the covariance of $\w\sim\prob[\w]{}$. Thus,
    \begin{gather*}
        C_{t+1} = \x_{t+1}^\intercal Q_K^c \x_{t+1} -\x_t^\intercal A_K^\intercal Q_K^c A_K\x_t 
        - \tr{Q_K^c H \Sigma_W H^\intercal}
    \end{gather*}
    Now recall $S_{t} \coloneqq \sum_{s=1}^{t} C_s$ and therefore, the cyclic property of trace and definition of $\Gamma_t$ imply that
    \begin{multline}\label{eq:S-t-expression}
        S_{t} = \tr{Q_K^c \Gamma_{t}} 
        -\tr{A_K^\intercal Q_K^cA_K(\Gamma_{t-1} + x_0 x_0^\intercal)}  \\
        - t\, \tr{Q_K^c H \Sigma_W H^\intercal}
    \end{multline}
    Also, by Lyapunov equation \cref{eq:Sigma-K}, we have the identity 
    \[\tr{Q_K^c\Sigma_K -A_K^\intercal Q_K^c A_K \Sigma_K - Q_K^c H \Sigma_W H^\intercal} = 0.\]
    Therefore, by applying the LLN in \cite[Corollary 12]{talebi_uniform_2024} to \cref{eq:S-t-expression} we conclude that as $t\to\infty$, 
        \( S_t/t \xrightarrow{a.s.}  0.\)
    Next, by applying the same Lyapunov identity, we can rewrite \cref{eq:S-t-expression} as 
    \begin{multline*}
        S_{t} = \tr{(Q_K^c -A_K^\intercal Q_K^cA_K) (\Gamma_{t}-t\Sigma_K)}  \\
        +\tr{A_K^\intercal Q_K^cA_K(\x_t \x_t ^\intercal - x_0 x_0^\intercal)}.
    \end{multline*}
    Recall that the CLT in \cite[Corollary 12]{talebi_uniform_2024} implies the convergence in distribution of $\frac{1}{\sqrt{t}}(\tr{M(\Gamma_{t} - t\Sigma_K)}$ for any constant matrix $M$. Furthermore, \Cref{lem:ave-conv} implies that $\tr{A_K^\intercal Q_K^c A_K \x_t \x_t^\intercal}/\sqrt{t}$ converges to zero in probability as $t\to\infty$.
    The first claim then follows by considering the linear (and thus continuous) mapping $\Gamma \mapsto \tr{(Q_K^c - A_K^\intercal Q_K^c A_K)\Gamma}$ and applying Continuous Mapping Theorem to this expression of $S_t$.
    Finally, we show the convergence of $N_t/t$ in \cref{eq:cond-var-def}.
    For that, $C_{t+1}$ is rewritten as 
    \( C_{t+1} = 2 (A_K \x_t)^\intercal Q_K^c H \w_{t+1} 
        + \tr{Q_K^c H (\w_{t+1} \w_{t+1}^\intercal - \Sigma_W )H^\intercal},
    \)
    and thus
    \begin{gather*}
    \begin{aligned}
        &\E[C_{t+1}^2|\F_{t}] = 4 (A_K \x_t)^\intercal Q_K^c H \Sigma_W H^\intercal Q_K^c (A_K \x_t) \\
        &+ 4 \x_t^\intercal A_K^\intercal Q_K^c
        \E[H\w_{t+1}\tr{Q_K^c H (\w_{t+1} \w_{t+1}^\intercal - \Sigma_W )H^\intercal} |\F_{t}] \\
        &+\E[\tr{Q_K^c H (\w_{t+1} \w_{t+1}^\intercal - \Sigma_W )H^\intercal}^2|\F_{t}]\\
        =& 4 \tr{ A_K^\intercal Q_K^c H \Sigma_W H^\intercal Q_K^c A_K \x_t \x_t^\intercal}  + 4 M_3^\intercal Q_K^c A_K \x_t  + m_4
    \end{aligned}
    \end{gather*}
    where we dropped the conditionals because $\w_{t+1}$ is independent of $\F_t$, with $m_4 = m_4[Q_K^c]$ as in the statement and \( M_3 = M_3[Q_K^c] \coloneq \E\big[H\w_{1}\tr{Q_K^c H (\w_{1} \w_{1}^\intercal - \Sigma_W )H^\intercal}\big]\), which are well-defined (bounded) by the moment condition on the noise.  
    %
    %
    Therefore, by definition of $N_t$ and \cite[Corollary 12]{talebi_uniform_2024}, we obtain the almost sure convergence: 
    \[
        \frac{1}{t}N_t \xrightarrow{a.s.} 4\tr{A_K^\intercal Q_K^c H \Sigma_W H^\intercal Q_K^c A_K \Sigma_K} + m_4[Q_K^c].
    \]
    But, by cyclic permutation property of the trace and the Lyapunov equation:
    \(
        \tr{A_K^\intercal Q_K^c H \Sigma_W H^\intercal Q_K^c A_K \Sigma_K} 
        = \tr{Q_K^c H \Sigma_W H^\intercal Q_K^c (\Sigma_K-H \Sigma_W H^\intercal)}.
    \)
    Combining the last two equations completes the proof.
\end{proof}

\section{Quadratic Ergodic-Risk COCP}\label{sec:primal-dual}
Herein, we show how the ergodic-risk criteria can be incorporated as a constraint in the optimal control framework posed in \cref{eq:optimization}.
The policy optimization (PO) approach to control design pivots on a parameterization of the feasible policies for the synthesis problem. 
One can view the \ac{lqr} cost naturally as a map $J(K) : K \mapsto J(\u=K \x)$. 
Also, for any stabilizing policy $K \in \mathcal{S}$, by using cyclic permutation property of trace, together with definitions in \cref{eq:def-s-t-Gamma-t}, we can compute the cost as $J_T(K) =\tr{Q_K\left(\x_0 \x_0^\intercal + \Gamma_T\right)}$
where $Q_K \coloneqq Q+K^\intercal R K$. Now, by \Cref{lem:ave-conv} and \Cref{thm:C-infty-N-convergence} we obtain that
\[J(K) = \lim_{T\to \infty} E[J_T(K)/T] = \tr{Q_K \Sigma_K},\]
with $\Sigma_K$ in \cref{eq:Sigma-K}. Note that $J(K)$ does not depend on the distribution of $\x_0$--as long as it has bounded second moment.

Recall, \Cref{thm:C-infty-N-convergence} ensures that $C_\infty$ is indeed distributed normally as long as $K$ is stabilizing, $(A_K,H)$ is controllable, and the noise has finite fourth moment; in this case $C_\infty \sim \mathcal{N}(0,\gamma_M^2)$.
But, any reasonable choice of a coherent risk measure on $C_\infty$ (such as Conditional Value-at-Risk CVaR$_{\alpha}$ and the Entropic Value-at-Risk EVaR$_{\alpha}$ on $C_\infty$ at a level $\alpha$), essentially reduces to an upperbound on a linear functional of $\gamma_M^2$. As discussed in \Cref{rem:asymp-var-conditional-var}, the asymptotic conditional variance $\gamma_N^2(K)$ can be interpreted as an ``estimate'' of $\gamma_N^2(K)$ which has a more tractable expression. Thus, herein we only consider constraints on $\gamma_N^2(K)$ and defer the other one to the extended version.
Therefore, the problem in \cref{eq:optimization} with the constraints on $\gamma_N^2$ over linear policies reduces to
\begin{align}\label{eq:optimization-reform}
    &\displaystyle \min \; J(K)  = \tr{Q_K \Sigma_K} \\
   &\text{s.t.} ~\gamma_N^2(K) \leq \bar\beta,\quad  
   K \in \stableK \cap \{K: (A_K,H) \text{ controllable}\}, \nonumber
\end{align}
with $\Sigma_K$ solving \cref{eq:Sigma-K}, $\gamma_N^2(K)$ defined in \Cref{thm:C-infty-N-convergence}, and for a given constant $\bar\beta$ encapsulating the risk level.

In this rest of this section, we develop a primal-dual algorithm to solve \cref{eq:optimization-reform} using Lagrange duality. Let us consider the Lagrangian $L:\mathcal{S}  \times \R_{\geq 0} \mapsto \R$ defined as
\begin{align}
    L(K, \lambda) &\coloneqq \tr{Q_K \Sigma_K} + \lambda \left(\gamma_N^2(K)-\bar\beta  \right) \label{eq:lagrangian}\\
     &= \tr{(Q_K + 4\lambda Q_K^c H \Sigma_W H^\intercal Q_K^c) \Sigma_K} + \lambda \beta[Q_K^c], \nonumber
\end{align}
with 
\(
    \beta[Q_K^c] \coloneqq -4\tr{(Q_K^c H \Sigma_W H^\intercal)^2 } + m_4[Q_K^c] -\bar\beta,
\)
and $m_4[Q_K^c]$ in \Cref{thm:C-infty-N-convergence}. Hereafter, we focus on the case where the risk functional $g$ does not depend on the input explicitly\footnote{This relates with the setting studied in \cite{tsiamis_risk-constrained_2020, zhao_global_2023}, where, in addition, the assumptions $Q^c = Q$ and $H = I$ are imposed.}.

\subsection{Quadratic Ergodic-Risk Criteria with $R^c=0$} \label{subseq:R-c-0}
Let us assume that the risk measure $g$ does not depend on control input explicitly; i.e. $R^c=0$, and so $Q_K^c = Q^c$ is constant. 
In addition to system-theoretic assumptions in \Cref{assmp:stability} that is necessary for feasibility of the optimization (non-empty domain $\stableK$), in the following result we also need controllability of the pair $(A^\intercal, Q^{\frac{1}{2}})$ to guarantee regularity of the Lagrangian, i.e. coercivity of $L(\cdot,\lambda)$ for each $\lambda$.
For simplicity, we consider slightly stronger conditions:

\begin{assumption}\label{asmp:Q-H-full-row-rank}
    Assume $Q \succ 0 $ and $H$ is full row rank.
\end{assumption}

The next result is obtain by leveraging  \cite[Lemma 1 and 2]{talebi_data-driven_2023}, where a detailed proof is provided in \cite{talebi_uniform_2024}.

\begin{lemma}\label{lem:Lagrangian}
Suppose Assumptions \ref{assmp:noise}, \ref{assmp:stability}, and \ref{asmp:Q-H-full-row-rank} hold.
For each $\lambda\geq0$, consider the Lagrangian $L_\lambda(\cdot) = L(\cdot,\lambda):\stableK \to \R$.
The following statements are true:
\begin{inparaenum}[(i)]
    \item $L_\lambda(\cdot)$ and $\gamma_N^2(\cdot)$ are smooth with 
    \(
    \nabla L_\lambda(K) = 2(RK + B^\intercal P_{(K,\lambda)} A_K)\Sigma_K 
    \)
    where $P_{(K,\lambda)}$ is the unique solution to
    \[P_{(K,\lambda)} = A_K^\intercal P_{(K,\lambda)} A_K + Q_K + 4\lambda Q^c H \Sigma_W H^\intercal Q^c.\]
    \item $L_\lambda(\cdot)$ is coercive with compact sublevel sets $\stableK_\alpha$ for each $\alpha>0$.
    \item $L_\lambda(\cdot)$ admits a unique global minimizer $K^*(\lambda) = \arg \min_{K\in\stableK} L_\lambda(K)$ that is stabilizing, given by
    \[K^*(\lambda) = -(R + B^\intercal P_{(K^*(\lambda),\lambda)}B)^{-1}B^\intercal P_{(K^*(\lambda),\lambda)} A,\]
    and $L_\lambda(K^*(\lambda)) = \tr{P_{(K^*(\lambda),\lambda)} H \Sigma_W H^\intercal} + \lambda \beta[Q^c]$.
    \item The restriction $L_\lambda(\cdot) |_{\stableK_\alpha}$ for any (non-empty) sublevel set $\stableK_\alpha$ has Lipschitz continuous gradient, is gradient dominated, and has a quadratic lower model; in particular, for all $K,K'\in\stableK_\alpha$ the following inequalities are true: 
    \(\|\nabla L_\lambda(K) - \nabla L_\lambda(K')\|_F \leq \ell \|K - K'\|_F\)
    and
    \(c_2 \|K-K^*\|_F^2\leq c_1[L_\lambda(K) - L_\lambda(K^*)] \leq \|\nabla L_\lambda(K)\|_F^2,\)
    for some positive constants $\ell(\alpha),c_1(\alpha)$, and $c_2(\alpha)$ that only depend on $\alpha$ and are independent of $K$.
\end{inparaenum}

\end{lemma}


\subsection{Strong Duality}\label{sec:strong-duality}
Now that $K^*(\lambda)$ is uniquely well-defined for each $\lambda \geq 0$, the \emph{dual problem} can be written in following forms
\begin{multline*}
    \sup_{\lambda\geq0} \min_{K \in \stableK} L(K,\lambda) 
= \sup_{\lambda\geq0} \tr{P_{(K^*(\lambda),\lambda)} H \Sigma_W H^\intercal} + \lambda \beta[Q^c]\\
=  \sup_{\lambda\geq0} \tr{Q_{K^*(\lambda)} \Sigma_{K^*(\lambda)}} + \lambda \left(\gamma_N^2(K^*(\lambda))-\bar\beta  \right),
\end{multline*}
with $P_{(K^*(\lambda),\lambda)}$ defined in \Cref{lem:Lagrangian}.
It is standard to assume the primal problem is strictly feasible: 
\begin{assumption}[Slater's Condition]\label{asmp:slater}
    Assume $\bar\beta$ is large enough such that there exists $K\in\stableK$ with $\gamma_N^2(K) < \bar\beta$.
\end{assumption}
This enables us to establish \emph{strong duality} for $L$ meaning that both primal and dual problems are feasible with identical values.
Because we know the cost function is globally lower bounded, $\tr{Q_{K} \Sigma_{K}} \geq 0$, then Slater's condition implies feasibility of the dual problem, and that there exists a finite $\lambda_0\geq 0$ such that $\gamma_N^2(K^*(\lambda_0)) \leq \bar\beta$. Now, if we let 
\begin{equation}\label{eq:lambda-star}
\lambda^* = \min\{\lambda_0 \;:\; \lambda_0 \geq 0 \text{ and } \gamma_N^2(K^*(\lambda_0)) \leq \bar\beta\}, 
\end{equation}
then we claim that the pair $(K^*(\lambda^*), \lambda^*)$ is the saddle point of the Lagrangian $L(K,\lambda)$, and therefore obtain the strong duality. For this claim to hold, by the Karush–Kuhn–Tucker (KKT) conditions, it suffice to show that:
\begin{equation}\label{eq:comp-slack}
    \mathrm{CS}(\lambda^*) \coloneqq \lambda^*(\gamma_N^2(K^*(\lambda^*)) - \bar\beta) = 0,
\end{equation}
aka the complementary slackness holds.
But, we know that $P_{(K,\lambda)}$ is real analytic in $(K,\lambda)$ and is positive definite for each $\lambda\geq0$. Thus, $K^*(\lambda)$ as defined in \Cref{lem:Lagrangian} is smooth in $\lambda$. So, $\gamma_N^2 \circ K^*(\cdot)$ is continuous by composition and lower bounded by zero. Therefore, complementary slackness follows because the minimum of a strictly monotone function on a compact set in \cref{eq:lambda-star} is attained at the boundary.

\subsection{The Primal-Dual Algorithm}\label{sec:primal-dual-algo}
By establishing a strong duality, we can solve the dual problem without loss of optimality. In particular, using the properties of $L_\lambda(\cdot)$ obtained in \Cref{lem:Lagrangian}, we devise simple primal-dual updates to solve \cref{eq:optimization-reform} by accessing the gradient and constrained violation values, as proposed in \Cref{algo}.
\begin{algorithm}[ht]
\caption{Primal-Dual Ergodic-Risk Constrained LQR}
    \begin{algorithmic}[1]
\State Set $K_0 \in \stableK_\alpha$ for some $\alpha>0$, $\lambda_0 =1$, and tolerance $\epsilon>0$ and stepsize $\eta_m =(\gamma_N^2(K_0) - \bar\beta)^{-1} (m+1)^{-1/2}$
\For{$m=0,1,\cdots, T=\mathcal{O}(\ln(\ln(\epsilon))/\epsilon^2)$}
\While{$\|\nabla_K L_{\lambda_m}(K)\|_F \geq \sqrt{\epsilon}$}
    \State $G \gets -(R + B^\intercal P_{K,{\lambda_m}} B)^{-1} \nabla L_{\lambda_m}(K) \Sigma_K^{-1}$
    \State $K \gets K + \frac{1}{2} G$ 
\EndWhile
\State $\lambda_{m+1} \gets \max\left[0, \lambda_{m} + \eta_m (\gamma_N^2(K) - \bar\beta) \right]$
\EndFor
\State \Return $(K,\lambda \coloneq \frac{1}{T}\sum_{m=1}^{T} \lambda_m)$
\end{algorithmic}
\label{algo}
\end{algorithm}
We next provide a convergence guarantee by combining recent LQR policy optimization \cite{fazel_global_2018, talebi_policy_2023} with saddle point optimization techniques \cite{nedic_subgradient_2009}, and illustrate the performance of \Cref{algo} through simulations.


\subsection{Convergence Guarantee}\label{subsec:convergence}

The results in \Cref{lem:Lagrangian} and the initialization in \Cref{algo} ensure that the premise of \cite[Theorem 4.3]{talebi_policy_2023} is satisfied. Also, $G$ is the Riemannian quasi-Newton direction and the update on $K$ is known as the Hewer's algorithm which is proved to converge at a quadratic rate, as discussed in detail in \cite[Remark 5]{talebi_policy_2023}. And therefore, the inner loop terminates very fast in $\ln(\ln(\epsilon))$ steps and essentially returns an $\epsilon$ accurate estimate of $K^*(\lambda_m)$. We could instead use a pure gradient descent algorithm, but this would result in a slower convergence rate.


Furthermore, the outer loop is expected to take $\mathcal{O}(1/\epsilon^2)$ steps to obtain $\epsilon$ error on the functional $L(K^*(\lambda), \lambda)$ following standard primal-dual guarantees \cite{nedic_subgradient_2009}. Thus, assuming that Assumptions \ref{assmp:noise}, \ref{assmp:stability}, \ref{asmp:Q-H-full-row-rank}, and \ref{asmp:slater} hold, \Cref{algo} obtains an $\epsilon$ accurate solution to problem \cref{eq:optimization-reform} in $\mathcal{O}(\ln(\ln(\epsilon))/\epsilon^2)$ steps.

\subsection{Simulations}\label{subsec:simulations}

Next, we compare the behavior of the ergodic-risk optimal policy with that of the LQR optimal policy on the Grumman X-29 aircraft dynamics, as studied in \cite{talebi_regularizability_2022}. We consider its longitudinal and lateral-directional dynamics in the Normal Digital Powered-Approach (ND-PA) mode, with a fixed discretization step size of $0.05$, following \cite[Tables 9 and 10]{bosworth_linearized_1992}. The system comprises four longitudinal states, four lateral-directional states, and five control inputs. Given the normalized state representation, we set $Q = Q_c = I_8$ and $R = R_c = I_5$.

The heavy-tailed process noise is drawn from a Student’s $t$-distribution with parameter $\nu = 5$, which has a finite fourth moment but unbounded fifth (and higher) moment. Additionally, to simulate external disturbances primarily affecting the unstable longitudinal dynamics, we introduce a longitudinal gust disturbance of magnitude $20$ at every $500$ time steps. As a consequence, standard risk-sensitive control methods such as Linear Exponential Quadratic Gaussian (LEQG), which rely on exponentiation of the cost functional, are inapplicable due to the non-existence of the required higher-order moments.

\begin{figure}[pt] \centering
\includegraphics[width=0.48\textwidth]{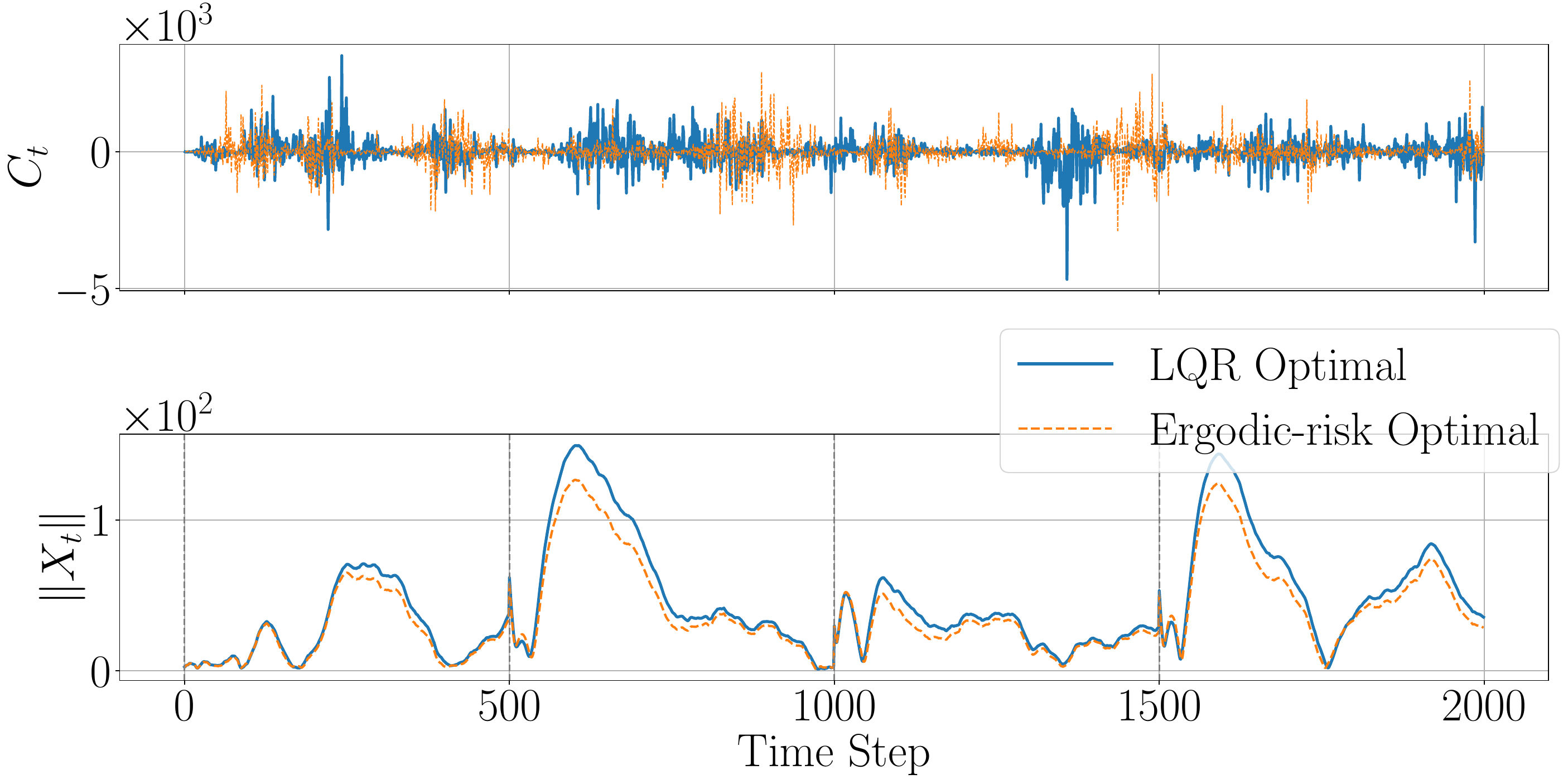} \caption{ \small Comparison of the optimal Ergodic-risk and optimal LQR policies for the Grumman X-29 aircraft under Student’s $t$-noise and simulated gust disturbances occurring every 500 time steps.} \label{fig:ergodic-vs-lqr-rollout} \vspace{-0.4cm} 
\end{figure}

The optimal ergodic-risk policy $K^*$, corresponding to a ergodic-risk level of $\bar{\beta} = 0.8 \gamma_N^2(K_{\text{LQR}})$, is computed using Algorithm~\ref{algo}. The resulting cost values are given by $J(K^*) = 623432$ versus $J(K_{\text{LQR}}) = 621829$. While the ergodic-risk policy \( K^* \) increases risk sensitivity by \( 20\% \) (quantified via the asymptotic conditional variance \( \gamma_N^2 \)), its average cost is only \( 0.25\% \) higher than that of the optimal LQR policy.
Also, a comparative evaluation of the ergodic-risk and LQR optimal policies is presented in \Cref{fig:ergodic-vs-lqr-rollout}. As illustrated, the ergodic-risk policy demonstrates superior resilience against gust disturbances relative to the LQR policy. The simulation codes are available at \cite{talebi_ergodic-risk_2025}.

\section{Conclusions}
We introduced ergodic-risk criteria in COCP as a flexible framework to account for long-term cumulative uncertainties. By incorporating linear constraints on $\gamma_N^2$ and leveraging recent advancements in policy optimization, we proposed a primal-dual algorithm with proven convergence guarantees. Key future directions of this work include further considering to develop policy optimization algorithms for directly constrain $\gamma_M^2$ where the risk functional also directly depends on the input signal, and developing sample-based algorithms.

\bibliographystyle{ieeetr}
\bibliography{citations}

\begin{thebibliography}{10}

\bibitem{rockafellar_optimization_2000}
R.~T. Rockafellar and S.~Uryasev, ``Optimization of conditional
  value-at-risk,'' {\em The Journal of Risk}, vol.~2, no.~3, pp.~21--41, 2000.

\bibitem{majumdar_how_2020}
A.~Majumdar and M.~Pavone, ``How {Should} a {Robot} {Assess} {Risk}? {Towards}
  an {Axiomatic} {Theory} of {Risk} in {Robotics},'' in {\em Robotics
  {Research}}, pp.~75--84, 2020.

\bibitem{eichler_risks_2013}
H.-G. Eichler and B.~Bloechl-Daum, et.~al., ``The risks of risk aversion in
  drug regulation,'' {\em Nature Reviews Drug Discovery}, vol.~12,
  pp.~907--916, Dec. 2013.

\bibitem{zhang_policy_2021}
K.~Zhang, B.~Hu, and T.~Başar, ``Policy {Optimization} for {H}-2 {Linear}
  {Control} with {H}-infinity {Robustness} {Guarantee}: {Implicit}
  {Regularization} and {Global} {Convergence},'' {\em SIAM Journal on Control
  and Optimization}, vol.~59, pp.~4081--4109, Jan. 2021.

\bibitem{whittle_risk-sensitive_1981}
P.~Whittle, ``Risk-sensitive linear/quadratic/gaussian control,'' {\em Advances
  in Applied Probability}, vol.~13, no.~4, pp.~764--777, 1981.

\bibitem{borkar_risk-constrained_2014}
V.~Borkar and R.~Jain, ``Risk-{Constrained} {Markov} {Decision} {Processes},''
  {\em IEEE Trans. Autom. Control}, pp.~2574--2579, Sept. 2014.

\bibitem{chow_risk-constrained_2018}
Y.~Chow, M.~Ghavamzadeh, L.~Janson, and M.~Pavone, ``Risk-constrained
  reinforcement learning with percentile risk criteria,'' {\em Journal of
  Machine Learning Research}, vol.~18, no.~167, pp.~1--51, 2018.

\bibitem{sopasakis_risk-averse_2019}
P.~Sopasakis, D.~Herceg, A.~Bemporad, and P.~Patrinos, ``Risk-averse model
  predictive control,'' {\em Automatica}, pp.~281--288, Feb. 2019.

\bibitem{ruszczynski_risk-averse_2010}
A.~Ruszczyński, ``Risk-averse dynamic programming for {Markov} decision
  processes,'' {\em Math. Program.}, pp.~235--261, Oct. 2010.

\bibitem{ahmadi-javid_entropic_2012}
A.~Ahmadi-Javid, ``Entropic value-at-risk: {A} new coherent risk measure,''
  {\em J. Optim. Theory Appl.}, vol.~155, pp.~1105--1123, Dec. 2012.

\bibitem{biswas_ergodic_2023}
A.~Biswas and V.~S. Borkar, ``Ergodic risk-sensitive control—{A} survey,''
  {\em Annual Reviews in Control}, vol.~55, pp.~118--141, Jan. 2023.

\bibitem{kishida_risk-aware_2023}
M.~Kishida and A.~Cetinkaya, ``Risk-aware linear quadratic control using
  conditional value-at-risk,'' {\em IEEE Transactions on Automatic Control},
  vol.~68, pp.~416--423, Jan. 2023.

\bibitem{tsiamis_risk-constrained_2020}
A.~Tsiamis, D.~S. Kalogerias, L.~F. Chamon, A.~Ribeiro, and G.~J. Pappas,
  ``Risk-constrained linear-quadratic regulators,'' in {\em 59th {IEEE}
  {Conference} on {Decision} and {Control}}, pp.~3040--3047, IEEE, 2020.

\bibitem{zhao_global_2023}
F.~Zhao, K.~You, and T.~Başar, ``Global convergence of policy gradient
  primal–dual methods for risk-constrained {LQRs},'' {\em IEEE Transactions
  on Automatic Control}, vol.~68, no.~5, pp.~2934--2949, 2023.

\bibitem{meyn_markov_2009}
S.~Meyn and R.~L. Tweedie, {\em Markov {Chains} and {Stochastic} {Stability}}.
\newblock Cambridge University Press, 2009.

\bibitem{durrett_probability_2019}
R.~Durrett, {\em Probability: {Theory} and {Examples}}.
\newblock Cambridge: Cambridge University Press, 5th~ed., 2019.

\bibitem{goodwin_control_2001}
G.~C. Goodwin, S.~F. Graebe, and M.~E. Salgado, {\em Control {System}
  {Design}}.
\newblock Prentice Hall, 2001.

\bibitem{talebi_uniform_2024}
S.~Talebi and N.~Li, ``Uniform {Ergodicity} and {Ergodic}-{Risk} {Constrained}
  {Policy} {Optimization},'' Sept. 2024.
\newblock arXiv:5856500.

\bibitem{komorowski_central_2012}
T.~Komorowski and A.~Walczuk, ``Central limit theorem for {Markov} processes
  with spectral gap in the {Wasserstein} metric,'' {\em Stochastic Processes
  and their Applications}, vol.~122, pp.~2155--2184, May 2012.

\bibitem{talebi_data-driven_2023}
S.~Talebi, A.~Taghvaei, and M.~Mesbahi, ``Data-driven optimal filtering for
  linear systems with unknown noise covariances,'' in {\em Advances in {Neural}
  {Inform}. {Process}. {Sys}.}, pp.~69546--69585, 2023.

\bibitem{hall_martingale_1980}
P.~Hall and C.~C. Heyde, {\em Martingale {Limit} {Theory} and {Its}
  {Application}}.
\newblock Academic Press, Dec. 1980.

\bibitem{fazel_global_2018}
M.~Fazel, R.~Ge, S.~Kakade, and M.~Mesbahi, ``Global convergence of policy
  gradient methods for the linear quadratic regulator,'' in {\em Int. {Conf}.
  on {Machine} {Learning}}, pp.~1467--1476, PMLR, July 2018.

\bibitem{talebi_policy_2023}
S.~Talebi and M.~Mesbahi, ``Policy optimization over submanifolds for linearly
  constrained feedback synthesis,'' {\em IEEE Transactions on Automatic
  Control}, pp.~1--16, 2023.

\bibitem{nedic_subgradient_2009}
A.~Nedić and A.~Ozdaglar, ``Subgradient methods for saddle-point problems,''
  {\em J. Optim. Theory Appl.}, vol.~142, pp.~205--228, July 2009.

\bibitem{talebi_regularizability_2022}
S.~Talebi, S.~Alemzadeh, N.~Rahimi, and M.~Mesbahi, ``On {Regularizability} and
  its application to online control of unstable {LTI} systems,'' {\em IEEE
  Trans on Automatic Control}, vol.~67, no.~12, pp.~6413--6428, 2022.

\bibitem{bosworth_linearized_1992}
J.~T. Bosworth, {\em Linearized aerodynamic and control law models of the
  {X}-{29A} airplane and comparison with flight data}, vol.~4356 of {\em {NASA}
  {Technical} {Memorandum}}.
\newblock National Aeronautics and Space Administration (NASA, Office of
  Management), 1992.

\bibitem{talebi_ergodic-risk_2025}
S.~Talebi, ``Ergodic-risk {Policy} {Optimization},'' Feb. 2025.
\newblock https://github.com/shahriarta/Ergodic-risk-Policy-Optimization.

\end{thebibliography}

\end{document}